\documentstyle[12pt]{article}
\def\s{\sigma} \def\etal{{\em et
al.\ }}\def\ni{\noindent}\def\ch{\choose}
\def\es{\emptyset}\def\sm{\setminus}\def\subq{\subseteq}\def\sub{\subset}
\def\s{\sigma}\def\Da{\Delta}
\def\cl{\centerline} \def\F{{\cal F}}\def\C{{\cal C}}\def\D{{\cal
D}}
 
\def\hD{\hat{\cal D}}\def\P{{\cal P}}\def\X{{\cal X}}
\def\sv{\s^v_{p,q}}\def\ss10{\s^v_{1,0}}\def\sso{\s^v_{1,q}}\def\sst{\s^v_{2,q}}

\def\lf{\lfloor}\def\rf{\rfloor}\def\lc{\lceil}\def\rc{\rceil}
 
\newtheorem{prop}{Proposition}\newtheorem{theo}{Theorem}\newtheorem{lem}{Lemma}
\newtheorem{defin}{Definition}\newtheorem{cor}{Corollary}
 
\begin{document} \title{Optimal
Pooling Designs with Error Detection} \author{David J. Balding
\\ School of Mathematical Sciences \\ Queen Mary \& Westfield
College\\ University of London\\ Mile End Road, London E1 4NS,
UK\vspace{5mm}\\David C. Torney \\ Theoretical Biology and Biophysics
\\ T-10 Mail Stop K-710 \\ Los Alamos National Laboratory \\Los Alamos
NM 87545 USA} \maketitle\thispagestyle{empty}
 
\cl{\it Abstract}\vspace{5mm}
 
{\narrower Consider a collection of objects, some of which may be
`bad', and a test which determines whether or not a given
sub-collection contains no bad objects. The non-adaptive pooling (or
group testing) problem involves identifying the bad objects using the
least number of tests applied in parallel. The `hypergeometric' case
occurs when an upper bound on the number of bad objects is known {\em a
priori}. Here, practical considerations lead us to impose the
additional requirement of {\em a posteriori} confirmation that the
bound is satisfied. A generalization of the problem in which occasional
errors in the test outcomes can occur is also considered.  Optimal
solutions to the general problem are shown to be equivalent to
maximum-size collections of subsets of a finite set satisfying a union
condition which generalizes that considered by Erd\"os \etal
\cite{erd}. Lower bounds on the number of tests required are derived
when the number of bad objects is believed to be either 1 or 2.
Steiner systems are shown to be optimal solutions in some
cases.\smallskip}\newpage
 
\section{Introduction}
 
Each of $n$ objects has an unknown binary status, `good' or `bad'. A
test is available which, except for occasional failures or errors,
establishes whether or not all the objects in a given collection are
good. The problem is to resolve the status of each object using the
minimum number of tests applied in parallel. The corresponding adaptive
problem, in which the choice of test at any stage can depend on the
outcomes of previous tests, is sometimes known as `group testing' (Wolf
\cite{wol}).
 
The objects may, for example, be electronic devices which can be tested
in series. Another example involves items in a database which are
categorized by a sequence of binary classifications and the task is to
partition the objects according to the $i$th classification. The
problem is formally similar to that of devising optimal
error-correcting codes using parity checks, except that here the test
result is `at least one bad object' rather than `an odd number of bad
objects'. Our work is motivated by an optimal design problem for
large-scale experiments aimed at constructing physical maps of human
chromosomes: the objects are chromosome fragments which are `bad' if
they contain a certain DNA sequence. An experimental test known as the
Polymerase Chain Reaction can determine whether or not a collection of
chromosome fragments are all good. In order to facilitate automation,
it is desirable that the experiments be applied in parallel so that the
experimental design is non-adaptive, or one-stage. Here, we derive
experimental designs which, with high probability, are one-stage
solutions to an appropriate formalization of the problem. These designs
may form stages in solutions to more general problems, for example
adaptive (multi-stage) designs which are optimal subject to a cost
function which penalizes additional stages.
 
A {\em pool} is a set of objects and a {\em design} is a set of pools.
Given a design $\D$, let $v$ denote the number of pools, so that
$v\equiv|\D|$. We will say that a pool is {\em good} if all the objects
in it are good, otherwise it is {\em bad}. Let $P$ denote the total
number of bad objects. The test usually distinguishes good pools from
bad, but we will also allow the possibility that for some pools the
test fails to produce a result and write $Q$ for the number of pools in
$\D$ which fail. Before applying the tests $P$ and $Q$ are unknown, but
we may have some prior information about them. One simple design
consists of testing each object individually a fixed number of times.
However if both $P\ll n$ and $Q\ll v$ then `better' designs are
possible.
 
There are several reasonable optimality criteria for $\D$. An
appropriate choice will depend in part on the prior knowledge of $P$
and $Q$. Bush \etal \cite{bus} and Hwang \& S\'os \cite{hwa} discuss
non-adaptive group testing in the `hypergeometric' case, in which $Q=0$
and $P$ is bounded above by a known constant $p$. They define $\D$ to
be an optimal solution if it maximizes $n$ for fixed $v$ among designs
such that the status of each object can be inferred from the pool
outcomes. The hypergeometric formulation has the drawback that it
assumes that the event $P>p$ is excluded {\em a priori}. It is not in
general possible to confirm {\em a posteriori} that $P\le p$ and hence
false conclusions may be drawn if, unexpectedly, $P>p$. In practice, a
large value of $p$ must be chosen to exclude this possibility. Here, we
modify the hypergeometric case by imposing the additional requirement
that the event $P>p$ can be distinguished {\em a posteriori}.
Consequently, it will be reasonable in practice to allow a small prior
probability that $P>p$.  Typically, lower values of $p$ can be chosen
than under the hypergeometric formulation and hence more efficient
designs constructed. The price for these advantages is that the designs
are not strictly non-adaptive: with small probability a second stage
will be required.
 
Allowing also for up to $q$ failures, we define $\D$ to be an optimal
solution if it maximizes $n$ for fixed $v$ subject to the requirement
that whenever $Q\le q$ we can infer from the pool outcomes either the
status of each object or that $P>p$. Proposition 1 of section
\ref{stat} establishes that optimal solutions $\D$ are equivalent to
maximum-size collections of subsets of a $v$-set such that every subset
in the collection has more than $q$ elements distinct from any union of
up to $p$ others. This condition is equivalent to that of $q$-error
detection and hence optimal $q$-failure designs are also optimal
$q$-error-detecting designs. In the $q=0$ case, we require that no
subset in the collection is contained in the union of $p$ others.
Hwang \& S\'os \cite{hwa} showed that this requirement characterizes
the $p$-complete designs defined by Bush \etal \cite{bus}.
 
In Theorems 1 and 2 we establish lower bounds on $v$ as a function of
$n$ for $p=1$ and $2$ and all $q\ge0$. The bounds coincide in some
cases with the sizes of certain Steiner system solutions which hence
are optimal. These results extend the results of Erd\"os \etal
\cite{erd} who considered the case $p=2$ and $q=0$. These authors
initially constrained the designs to be uniform, that is each object
occurs in the same number of pools. They subsequently derived an
asymptotic bound in the unconstrained case.  Here, we do not require
uniformity but we note that the bounds given in Theorems 1 and 2 can
only be achieved by uniform designs.  Ruszink\'o \cite{rus} derives
asymptotic bounds for $q=0$ and arbitrary $p$, but in the case $p=2$
the bound obtained by Erd\"os \etal \cite{erd} is tighter.
 
\section{Definitions and statement of results}\label{stat} For positive
integers $0\le i\le j$, let $\X_j$ denote the set of subsets of
$\{1,2,\ldots,j\}$ and define
\begin{equation}\X^i_j\equiv\{C\in\X_j:|C|=i\}.\end{equation} Pools are
elements of $\X_n$ and designs are subsets of $\X_n$. Given a design
$\D\equiv\{A_1,\ldots,A_v\}$ we will write
$\hD\equiv\{B_1,\ldots,B_n\}$ for the {\em dual} of $\D$ defined by
$i\in A_j$ if and only if $j\in B_i$. Thus $A_j$ indexes the objects in
the $j$th pool whereas $B_i$ indexes the pools which contain the $i$th
object. Let $\phi(A)$ denote the set of indices of bad pools in $\D$
when the objects indexed by $A$ are bad and no failures occur, that is
\begin{equation}\phi(A)=\cup_{i\in A}B_i.\end{equation} We say that
$\D$ is a $p$-bad, $0$-failure solution, or $(p,0)$-solution, if from
$\phi(A)$ we can infer either $A$ or that $|A|>p$, assuming that no
failures occur. This occurs if and only if
\begin{equation}\phi(A)\ne\phi(A')\quad\hbox{for all }A,A'\in\X_n
\hbox{ such that }A\ne A'\hbox{ and }|A|\le
p.\label{def0}\end{equation} Note that in the hypergeometric case
(Hwang \& S\'os, \cite{hwa}), $\phi(A)\ne\phi(A')$ is required only
when both $|A|\le p$ and $|A'|\le p$.
 
We define $\D$ to be a $p$-bad, $q$-failure solution, or
$(p,q)$-solution, if from $\phi(A)$ we can infer either $A$ or that
$|A|>p$, even in the presence of up to $q$ failures. This occurs if and
only if each $(v{-}q)$-subset of $\D$ is a $(p,0)$-solution. We write
$\sv$ for the set of duals of $(p,q)$-solutions and say that $\D$ is
optimal if $\hD$ has maximum cardinality in $\sv$. From (\ref{def0}) it
follows that $\D$ is a $(p,q)$-solution if and only if
\begin{equation}|\phi(A)\Da\phi(A')|>q\quad\hbox{for all }A,A'\in\X_n
\hbox{ such that }A\ne A'\hbox{ and }|A|\le
p,\label{defq}\end{equation} in which $B\Da
C\equiv(B{\sm}C)\cup(C{\sm}B)$. Note that (\ref{defq}) can be regarded
as the definition of a solution in the case that test failures do not
occur but up to $q$ wrong outcomes may be recorded and the detection of
any such error is required. Hence optimal $p$-bad, $q$-failure
solutions are also optimal $p$-bad, $q$-error-detecting solutions.
\begin{prop} A design $\D$ is a $(p,q)$-solution, that is $\hD\in\sv$,
if and only if \begin{equation}\bigl|B_i\sm\phi(A)\bigr|>q\quad\hbox{
for every } A\in\X_n\hbox{ with }|A|\le p\hbox{ and all }
i\in\{1,2,\ldots,n\}{\sm}A.\label{cond}\end{equation}\end{prop}
\begin{cor} A design $D$ satisfies $\hD\in\sso$ if and only if
$|B{\sm}B'|>q$ for all distinct $B,B'\in\hD$.\label{c1}\end{cor}
 
\ni{\em Proof}\quad By considering the case that $A\in\X^p_n$ and
$A'=A\cup\{i\}$ for some $i\notin A$, we see that (\ref{cond}) is
necessary for (\ref{defq}). Suppose that $A,A'\in\X_n$ with $A\ne A'$
and $|A|\le p$. If $A'{\sm}A\ne\es$ then it follows from (\ref{cond})
that $|\phi(A'){\sm}\phi(A)|>q$. Alternatively, if $A'{\sm}A=\es$ then
both $|A'|\le p$ and $A{\sm}A'\ne\es$ and hence (\ref{cond}) implies
that $|\phi(A){\sm}\phi(A')|>q$. In either case we have
$|\phi(A)\Da\phi(A')|>q$ and hence (\ref{cond}) is sufficient for
(\ref{defq}).\vskip3mm
 
Let $0\le t<k\le v$. A $(t,k,v)$-packing is a set $\P\subq\X^k_v$ such
that the intersection of any two elements of $\P$ has cardinality at
most $t$. A direct corollary of Proposition 1 is that if $\hD$ is a
$(t,pt{+}q{+}1,v)$-packing then $\hD\in\sv$. If $|\P|={v\ch t+1}{k\ch
t+1}^{-1}$ then each element of $\X^{t+1}_v$ is contained in precisely
one element of $\P$ and $\P$ is also called a Steiner system, denoted
$S(t{+}1,k,v)$.  For further details including a list of small Steiner
systems known to exist, we refer to Beth \etal \cite{bet}.
 
\begin{theo} If a design $\D$ satisfies $\hD\in\sso$ then
$n\equiv|\hD|$ satisfies \begin{equation}n\le{1\over K_q}{v\ch \lf
v/2\rf},\label{b1}\end{equation} in which $K_0=1$ and, for $q$ even,
\begin{equation}K_q=\sum_{s=0}^{q/2}{\lf v/2\rf\ch s}{\lc v/2\rc\ch s},
\label{kq1}\end{equation} while for $q$ odd,
\begin{equation}K_q=K_{q-1}+{1\over T}{\lf v/2\rf\ch(q{+}1)/2}{\lc
v/2\rc\ch(q{+}1)/2}\label{kq2}\end{equation} where $T\equiv\lf2\lf
v/2\rf/(q{+}1)\rf$.\end{theo} \begin{cor}[Sperner, 1928] The set
$\X^{\lf v/2\rf}_v$ is optimal in $\ss10$.\end{cor}\begin{cor}If $S(\lf
v/2\rf{-}1,\lf v/2\rf,v)$ exists then it is optimal in
$\s^v_{1,1}$.\end{cor}
 
\begin{theo} If a design $\D$ satisfies $\hD\in\sst$ then
\begin{equation}n\le{v\ch t^*}{2t^*+q-1\ch
t^*}^{-1},\label{b2}\end{equation} in which $t^*$ is the least integer
value of $t$ such that \begin{equation} v\le5t+2+{q(q{-}1)\over
t{+}q},\label{deft}\end{equation} so that $t^*=\lc(v{-}2)/5\rc$ if
$q=0$ or $1$.\end{theo} \begin{cor} If $S(t^*,2t^*{+}q{-}1,v)$ exists
then it is optimal in $\sst$.\end{cor}
 
In the $p=1$ case, Stirling's formula and (\ref{b1}) give
\begin{equation}n<\lf q/2\rf!\lc q/2\rc!{2^{v+q+1}\over v^q\sqrt{2\pi
v}}, \end{equation} and hence asymptotically
\begin{equation}v>{\log(n)\over\log(2)}\left(1{+}o(1)\right).\end{equation}
For $p=2$, Stirling's formula and (\ref{b2}) give
\begin{equation}n<{5^{v+1/2}\over2^{2v+q+1/2}},\end{equation} and hence
\begin{equation}v>{\log(n)\over\log(5/4)}\left(1{+}o(1)\right).
\label{ab2}\end{equation} The asymptotic bound (\ref{ab2}) was obtained
by Erd\"os \etal \cite{erd}.
 
\section{Proof of Theorem 1} By Corollary \ref{c1}, if $q=0$ then each
$\hD\in\ss10$ satisfies the requirement of no pairwise containment and
the theorem was proved in this case by Sperner \cite{spe}. A simple
proof of Sperner's result is given by Lubell \cite{lub}. Here, and in
the sequel, `chain' will always mean a maximal chain of $\X_v$, ordered
by inclusion. Each chain contains at most one element of $\hD$ and
hence we can associate with each $B\in\hD$ a `cost' which is the
proportion of all chains which contain $B$ and hence no further element
of $\hD$.  For any $k$, the set of chains can be partitioned by the
$k$-sets into $v\ch k$ equal parts. Therefore the cost of $B$ is
$1/{v\ch k}$, where $k\equiv|B|$, which is minimized at $k=\lf v/2\rf$.
Since $\X^{\lf v/2\rf}_v$ consists only of minimal cost elements but
achieves the maximum total cost of one, it is optimal in $\ss10$.
 
This argument can be extended to the $p=1$, $q>0$ case. Suppose first
that $q$ is even and consider $B\in\hD$ with $|B|=k$. Define the
$s$-neighbours of $B$ to be the sets $C\in\X_v^k$ such that
$|B{\sm}C|=|C{\sm}B|=s$. We will say that $B$ `blocks' the chains which
contain one of its $s$-neighbours for some $s\le q/2$. Let $B'$ denote
an element of $\hD$ distinct from $B$. If a chain contains both an
$s$-neighbour of $B$ and an $s'$-neighbour of $B'$ then either
$|B{\sm}B'|\le s{+}s'$ or $|B'{\sm}B|\le s{+}s'$. It follows from
Corollary \ref{c1} that we cannot have both $s\le q/2$ and $s'\le
q/2$.  Therefore a chain cannot be blocked by more than one element of
$\hD$.  Each chain blocked by $B$ can contain no element of $\hD$ other
than $B$ and hence we can associate with $B$ the cost $h(B)$ which is
the proportion of chains blocked by $B$. The value of $h(B)$ is
\begin{equation}h(B)=K_{q,k}{v\ch k}^{-1},\end{equation} where
$K_{q,k}$ denotes the total number of $s$-neighbours of $B$ with $s\le
q/2$, which is given by \begin{equation}K_{q,k}=\sum_{s=0}^{q/2}{k\ch
s}{v{-}k\ch s}.\end{equation} It is readily verified that $h(B)$ is
minimized when $k=\lf v/2\rf$ or $\lc v/2\rc$ and, in either case, the
value of $K_{q,k}$ is $K_q$, defined at (\ref{kq1}). Since chains
cannot be multiply blocked, the sum of the costs of the elements of
$\hD$ cannot exceed one. Therefore $n$ is bounded above by the inverse
of the minimal cost which establishes (\ref{b1}) in the case $q$ even.
 
Now suppose that $q$ is odd. Each chain which contains an $s$-neighbour
of $B$, for some $s\le(q{+}1)/2$, can contain no element of $\hD$ other
than $B$ and we say that these chains are blocked by $B$. If $B'$ is an
element of $\hD$ distinct from $B$ then no chain can be both an
$s$-neighbour of $B$ and an $s'$-neighbour of $B'$ when $s{+}s'\le q$.
However, if either $|B{\sm}B'|=q{+}1$ or $|B'{\sm}B|=q{+}1$ then there
exist chains which contain both a $((q{+}1)/2)$-neighbour of $B$ and a
$((q{+}1)/2)$-neighbour of $B'$ and hence such chains are multiply
blocked.
 
Let $\C$ denote the set of elements of $\hD$ which have a
$((q{+}1)/2)$-neighbour in the chain
$\{\emptyset,\{1\},\{1,2\},\ldots,\{1,2,\ldots,v\}\}$. If $B\in\C$ with
$|B|=k$ then
\begin{equation}|B\cap\{k{+}1,k{+}2,\ldots,v\}|={q+1\over2}=
|\{1,2,\ldots,k\}\sm B|.\end{equation} Further, if $B'\in\C$ and $B'\ne
B$ then, since both $|B{\sm}B'|>q$ and $|B'{\sm}B|>q$ must be
satisfied, $B'$ contains $\{1,2,\ldots,k\}\sm B$ and is disjoint from
$\{k{+}1,k{+}2,\ldots,v\}\cap B$. Hence $|\C|$ cannot exceed
min$\{\lf2k/(q{+}1)\rf,\lf2(v{-}k)/(q{+}1)\rf\}$, which takes maximum
value $T$ when $k=\lf v/2\rf$. Therefore the number of
$((q{+}1)/2)$-neighbours in $\hD$ of a given chain is at most $T$.
 
Associate with $B$ a cost \begin{equation}h'(B)=K'_{q,k}{v\ch
k}^{-1},\end{equation} where $K'_{q,k}$ is $K_{q-1,k}$ plus $1/T$ times
the number of $((q{+}1)/2)$-neighbours of $B$. This cost is minimized
at $k=\lf v/2\rf$ or $\lc v/2\rc$ and in either case $K'_{q,k}$ is
$K_q$ defined at (\ref{kq2}). The sum of the costs of the elements of
$\hD$ cannot exceed one and hence the theorem.
 
\section{Proof of Theorem 2} \begin{defin} For any $\hD\in\sst$, we
follow Erd\"os \etal {\em \cite{erd}} and say that $b\in\X_v$ is {\em
private in} $\hD$ if there exists a unique $B\in\hD$ such that $b\subq
B$.\label{d1}\end{defin}
 
\begin{defin} If $B\in\X_v$ with $|B|>q$ then $\F\sub\X_v$ is a {\em
(2,q)-cover} of $B$ precisely if both \begin{enumerate}\item if
$b\in\F$ and $b\sub b'\subq B$ then $b'\in\F$;\quad and \item for every
$b\subq B$ with $|B{\sm}b|\le q$ at least one part of every
two-partition of $b$ is in $\F$.\end{enumerate}\label{d2}\end{defin}
 
\begin{lem} A design $\D$ satisfies $\hD\in\sst$ if and only if for
each $B\in\hD$, the sets which are private in $\hD$ form a
$(2,q)$-cover of $B$.\label{l1}\end{lem}
 
\ni{\em Proof}\quad Let $b\subq B$ with $|B{\sm}b|\le q$. It follows
from Proposition 1 that $|B|>q$ and $b$ is private in $\hD$.  If there
exists a partition of $b$ into two non-private parts then there must be
$C$ and $C'$ in $\hD{\sm}\{B\}$ such that $b\subq(C\cup C')$ and hence
$|B\sm(C{\cup}C')|\le q$, which contradicts Proposition 1.  Therefore a
necessary condition for $\D\in\sst$ is that for each $B\in\hD$ the
private subsets of $B$ in $\hD$ form a $(2,q)$-cover.  Sufficiency is
immediate from Proposition 1.\vskip3mm
 
\begin{defin} For $\F$ a $(2,q)$-cover of $B\in\X_v$, define $h(B,\F)$
to be the proportion of all chains which intersect $\F$.\end{defin}
 
\begin{prop} For any $B\in\X_v$ with $|B|=k$ such that $q<k<v{-}1$ and
$\F$ a $(2,q)$-cover of $B$, \begin{equation}h(B,\F)\ge{2t{+}q{-}1\ch
t}{v\ch t}^{-1},\label{hb}\end{equation} in which
$t=\lf(k{-}q{+}1)/2\rf$.  Equality is achieved in (\ref{hb}) if and
only if $k{-}q$ is odd and $\F=\F^*$ where $\F^*\equiv\{b\sub B:|b|\ge
t\}$.\label{p2}\end{prop}
 
\ni{\em Proof}\quad If $k=2t{+}q{-}1$ then every two-partition of any
$(k{-}q)$-subset of $B$ contains a part $b$ such that $|b|\ge t$ and
hence $\F^*$ is a $(2,q)$-cover of $B$. Further, $h(B,\F^*)$ is
precisely the proportion of chains which contain a $t$-subset of $B$
and hence $h(B,\F^*)$ achieves equality in (\ref{hb}).
 
Suppose that $k=2t{+}q{-}1$ and let $s$ denote the largest integer such
there exists some $b\sub B$ with $|b|=t{+}s$ and $b\notin\F$. It
follows from Definition \ref{d2} that $s<t{-}1$. If $s<0$ then either
$\F=\F^*$ or $\F^*\sub\F$ and $h(B,\F^*)<h(B,\F)$ and hence we may
assume that $s\ge0$. If there exists some $b\in\F$ such that
$|b|<t{-}s{-}1$ then $\F{\sm}\{b\}$ is also a $(2,q)$-cover of $B$ and
$h(B,\F{\sm}\{b\})<h(B,\F)$. Hence we may assume that there is no such
$b$. From Definition \ref{d2}, if $b$ is a $(k{-}q)$-subset of $B$ then
the number of $(t{+}s)$-subsets of $b$ not in $\F$ is not greater than
the number of $(t{-}s{-}1)$-subsets of $b$ in $\F$.  Summing over all
such $b$, each $r$-set occurs in ${k{-}r\ch q}$ terms of the sum and
hence \begin{equation}{k{-}t{-}s\ch q}\bar f_{t+s}\le{k{-}t{+}s{+}1\ch
q} f_{t-s-1},\label{in}\end{equation} in which $\bar f_r$ denotes the
number of $r$-subsets of $B$ not in $\F$ while $f_r$ denotes the number
of such subsets in $\F$. Inequality (\ref{in}) is equivalent to
\begin{equation}\bar f_{t+s}\le
f_{t-s-1}{(k{-}t{+}s{+}1)!(t{-}s{-}1)!\over(k{-}t{-}s)!(t{+}s)!}.
\label{in2}\end{equation} Construct $\F'$ from $\F$ by removing all
$(t{-}s{-}1)$-sets and adding any missing $(t{+}s)$-subsets of $B$, so
that \begin{equation}\F'\equiv\{b\sub B:|b|=t{+}s\}\cup\F\sm\{b\sub
B:|b|=t{-}s{-}1\}.\end{equation} Now, $\F'$ is also a $(2,q)$-cover of
$B$ and $h(B,\F')-h(B,\F)$ is precisely the proportion of chains which
contain a $(t{+}s)$-subset of $B$ but no element of $\F$ minus the
proportion of chains which contain a $(t{-}s{-}1)$-set in $\F$ but no
other element of $\F$. Therefore
\begin{equation}h(B,\F')-h(B,\F)={v-k\over v{-}t{-}s}\bar f_{t+s}{v\ch
t{+}s}^{-1}-{v-k\over v{-}t{+}s{+}1}f_{t-s-1}{v\ch
t{-}s{-}1}^{-1},\end{equation} and hence $h(B,\F')\ge h(B,\F)$ if and
only if \begin{equation}\bar f_{t+s}\ge
f_{t-s-1}{(v{-}t{+}s)!(t{-}s{-}1)!
\over(v{-}t{-}s{-}1)!(t{+}s)!}.\label{con}\end{equation} Since
$k{+}1<v$, inequality (\ref{con}) contradicts (\ref{in2}) and hence
$h(B,\F')<h(B,\F)$. Therefore $h(B,\F)$ is not minimal and the
proposition is established in this case.
 
When $k=2t$ and $q=0$, for every partition of $B$ into two $t$-sets,
one of the parts is in $\F$. Hence $\F$ contains at least ${2t{-}1\ch
t}$ sets of size $t$. From an argument similar to that above, it is
readily shown that $h(B,\F)$ is minimized when $\F$ contains every
$(t{+}1)$-subset of $B$ but no $(t{-}1)$-subset. Hence the bound
(\ref{hb}) follows with strict inequality. If $k=2t{+}q$, $q>0$, then
for any $x\in B$ the set $\{b\in\F:x\notin b\}$ is a $(2,q{-}1)$-cover
of $B{\sm}\{x\}$ and the proposition follows from the case
$k=2t{+}q{-}1$.\vspace{3mm}
 
\ni{\em Proof of Theorem 2}\quad Let $x_t$ denote the RHS of
(\ref{hb}). Then \begin{equation}{x_{t+1}\over
x_t}={(2t{+}q{+}1)(2t{+}q)\over(v{-}t)(t{+}q)}
\label{maxt}\end{equation} which exceeds one if and only if inequality
(\ref{deft}) is not satisfied.  Hence $x_t$ is minimized at $t\equiv
t^*$. From Definition \ref{d1}, if $B\in\hD$ and $\hD\in\sst$ then any
chain which intersects $\F$ cannot intersect any other set which is
private in $\hD$. Therefore, invoking Lemma \ref{l1}, $n\equiv|\hD|$ is
bounded above by the inverse of the minimum value of $h(B,\F)$ over all
$\F$ a $(2,q)$-cover of $B$, which in turn is bounded above by
$1/x_{t^*}$.
 
\vspace{5mm}\ni{\bf Acknowledgements}\quad We thank Dr Charles Goldie
of QMW for helpful comments on an early draft of the manuscript. Work
supported in part by the UK Science and Engineering Research Council
under grants GR/F 98727 (DJB) and GR/J 05880 (DCT) and in part through
the Center for Human Genome Studies at Los Alamos National Laboratory
under grant US DOE/OHER ERWF118.
 
\end{document}